\def\R{\mathbb R}
\def\N{\mathbb N}
\def\Z{\mathbb Z}
\def\Q{\mathbb Q}
\def\pf{\begin{proof}}
\def\pfk{\end{proof}}

\documentclass[11pt]{article}
\usepackage{amssymb,amsmath,amsthm}

\newtheorem{thm}{Theorem}

\begin{document}




\title{\bf SELF-MATCHING PROPERTIES OF BEATTY SEQUENCES}


\author{Zuzana Mas\'akov\'a, Edita Pelantov\'a}
\date{}

\maketitle

\begin{center}
{Department of Mathematics, FNSPE, Czech Technical University\\
Trojanova 13, 120 00 Praha 2, Czech Republic\\
E-mail: masakova@km1.fjfi.cvut.cz, pelantova@km1.fjfi.cvut.cz}\\
\end{center}

\bigskip

\begin{abstract}
We study the selfmatching properties of Beatty sequences, in
particular of the graph of the function $\lfloor j\beta\rfloor $
against $j$ for every quadratic unit $\beta\in(0,1)$. We show that
translation in the argument by an element $G_i$ of generalized
Fibonacci sequence causes almost always the translation of the
value of function by $G_{i-1}$. More precisely, for fixed
$i\in\N$, we have $\bigl\lfloor \beta(j+G_i)\bigr\rfloor = \lfloor
\beta j\rfloor +G_{i-1}$, where $j\notin U_i$. We determine the
set $U_i$ of mismatches and show that it has a low frequency,
namely $\beta^i$.
\end{abstract}




{\section{Introduction}}

Sequences of the form $\bigl(\lfloor j\alpha\rfloor
\bigr)_{j\in\N}$ for $\alpha>1$, now known as Beatty sequences,
have been first studied in the context of the famous problem of
covering the set of positive integers by disjoint
sequences~\cite{beatty}. Further results in the direction of the
so-called disjoint covering systems are due
to~\cite{fraenkel,graham,tijdeman} and others. Other aspects of
Beatty sequences were then studied, such as their generation using
graphs~\cite{deBruijn}, their relation to generating
functions~\cite{Komatsu,OBryant}, their substitution
invariance~\cite{komatsu2,parvaix}, etc. A good source of
references on Beatty sequences and other related problems can be
found in~\cite{brown,stolarsky}.

In~\cite{bunder} the authors study the self-matching properties of
the Beatty sequence $\bigl(\lfloor j\tau\rfloor \bigr)_{j\in\N}$
for $\tau=\frac12(\sqrt5-1)$, the golden ratio. Their study is
rather technical; they have used for their proof the Zeckendorf
representation of integers as a sum of distinct Fibonacci numbers.
The authors also state an open question whether the results
obtained can be generalized to other irrationals than $\tau$. In
our paper we answer this question in the affirmative. We show that
Beatty sequences $\bigl(\lfloor j\alpha\rfloor \bigr)_{j\in\N}$
for quadratic Pisot units $\alpha$ have analogical self-matching
property, and for our proof we use a simpler method, based on the
cut-and-project scheme.

It is interesting to mention that Beatty sequences, Fibonacci
numbers and cut-and-project scheme attracted the attention of
physicists in recent years because of their applications for
mathematical description of non-crystallographic solids with
long-range order, the so-called quasicrystals, discovered in
1982~\cite{schecht}. The first observed quasicrystals revealed
crystallographically forbidden rotational symmetry of order 5.
This necessitates, for the algebraic description of the
mathematical model of such a structure, the use of the quadratic
field $\Q(\tau)$. Such a model is self-similar with the scaling
factor $\tau^{-1}$. Later, one observed existence of quasicrystals
with 8 and 12-fold rotational symmetries, corresponding to
mathematical models with selfsimilar factors $\mu^{-1}=1+\sqrt2$
and $\nu^{-1}=2+\sqrt3$. Note that all $\tau$, $\mu$, and $\nu$
are quadratic Pisot units, i.e.\ belong to the class of numbers
for which the result of Bunder and Tognetti is generalized here.

\bigskip
{\section{Quadratic Pisot units and cut-and-project scheme}}

The self-matching properties of the Beatty sequence $\bigl(\lfloor
j\tau\rfloor \bigr)_{j\in\N}$ are best displayed on the graph of
$\lfloor j\tau\rfloor $ against $j\in\N$. Important role is played
by the Fibonacci numbers,
$$
F_0=0,\ F_1=1,\quad F_{k+1}=F_k+F_{k-1},\ \hbox{for }k\geq 1.
$$
The result of~\cite{bunder} states that
\begin{equation}\label{eq:cislo}
\lfloor (j+F_i)\tau\rfloor =\lfloor j\tau\rfloor +F_{i-1}\,,
\end{equation}
except isolated mismatches of frequency $\tau^{i}$, namely at
points $j=kF_{i+1}+\lfloor k\tau\rfloor  F_i$.

Our aim is to show a very simple proof of the mentioned results
that is valid for all quadratic units $\beta\in(0,1)$. Every such
unit is a solution of the quadratic equation
$$
x^2+mx=1\,,\ m\in\N\,,\qquad\hbox{ or }\qquad x^2-mx=-1\,,\
m\in\N\,,\ m\geq 3\,.
$$
The considerations will slightly differ in the two cases.

\begin{itemize}
\item[(a)]
Let $\beta\in(0,1)$ satisfy $\beta^2+m\beta=1$ for $m\in\N$. The
algebraic conjugate of $\beta$, i.e.\ the other root of the
equation, satisfies $\beta'<-1$. We define the generalized
Fibonacci sequence
\begin{equation}\label{ega}
G_0=0\,,\quad G_1=1\,,\quad G_{n+2}=mG_{n+1}+G_n\,,\ n\geq 0\,.
\end{equation}
It is easy to show by induction that for $i\in\N$, we have
\begin{equation}\label{egaa}
(-1)^{i+1}\beta^i = G_i\beta-G_{i-1}\quad\hbox{ and }\quad
(-1)^{i+1}{\beta'}^i = G_i\beta'-G_{i-1}\,.
\end{equation}

\item[(b)]
Let $\beta\in(0,1)$ satisfy $\beta^2-m\beta=-1$ for $m\in\N$,
$m\geq 3$. The algebraic conjugate of $\beta$ satisfies
$\beta'>1$. We define
\begin{equation}\label{egb}
G_0=0\,,\quad G_1=1\,,\quad G_{n+2}=mG_{n+1}-G_n\,,\ n\geq 0\,.
\end{equation}
In this case, we have for $i\in\N$
\begin{equation}\label{egbb}
\beta^i = G_i\beta-G_{i-1}\quad\hbox{ and }\quad {\beta'}^i =
G_i\beta'-G_{i-1}\,.
\end{equation}
\end{itemize}

The proof we give here is based on the algebraic expression for
one-dimensional cut-and-project sets~\cite{c&p}. Let $V_1$, $V_2$
be straight lines in $\R^2$ determined by vectors $(\beta,-1)$ and
$(\beta',-1)$, respectively. The projection of the square lattice
$\Z^2$ on the line $V_1$ along the direction of $V_2$ is given by
$$
(a,b)=(a+b\beta')\vec{x}_1 + (a+b\beta)\vec{x}_2\,,\quad\hbox{ for
}(a,b)\in\Z^2\,,
$$
where $\vec{x}_1=\frac1{\beta-\beta'}(\beta,-1)$ and
$\vec{x}_2=\frac1{\beta'-\beta}(\beta',-1)$. For the description
of the projection of $\Z^2$ on $V_1$ it suffices to consider the
set
$$
\Z[\beta']:=\{a+b\beta'\mid a,b\in\Z\}\,.
$$
The integral basis of this free abelian group is $(1,\beta')$, and
thus every element $x$ of $\Z[\beta']$ has a unique expression in
this base. We will say that $a$ is the rational part of
$x=a+b\beta'$ and $b$ is its irrational part. Since $\beta'$ is a
quadratic unit, $\Z[\beta']$ is a ring and, moreover, it satisfies
\begin{equation}\label{e:unit}
\beta'\Z[\beta']=\Z[\beta']\,.
\end{equation}

A cut-and-project set is the set of projections of points of
$\Z^2$ to $V_1$, that are found in a strip of bounded width,
parallel to the straight line $V_1$. Formally, for a bounded
interval $\Omega$ we define
$$
\Sigma(\Omega) = \{ a+b\beta' \mid a,b\in\Z,\
a+b\beta\in\Omega\}\,.
$$
Note that $a+b\beta'$ corresponds to the projection of the point
$(a,b)$ to the straight line $V_1$ along $V_2$, whereas $a+b\beta$
corresponds to the projection of the same lattice point to $V_2$
along $V_1$.

Among the simple properties of cut-and-project sets that we use
here are
$$
\Sigma(\Omega-1) = -1+\Sigma(\Omega)\,,\qquad\qquad
\beta'\Sigma(\Omega) = \Sigma(\beta\Omega)\,,
$$
where the latter is a consequence of~\eqref{e:unit}. If the
interval $\Omega$ is of unit length, one can derive directly from
the definition a simpler expression for $\Sigma(\Omega)$. In
particular, we have
$$
\Sigma[0,1) = \bigl\{a+b\beta' \,\bigm|\, a+b\beta \in
[0,1)\bigr\} = \{b\beta'-\lfloor b\beta\rfloor  \mid b\in\Z\}\,,
$$
where we use that the condition $0\leq a+b\beta <1$ is satisfied
if and only if $a=\lceil -b\beta\rceil = -\lfloor b\beta\rfloor $.

Let us mention that the above properties of one-dimensional
cut-and-project sets, and many others, are explained in the review
article~\cite{c&p}.

\section{Self-matching property of the graph $\lfloor j\beta\rfloor $
against $j$} 

Important role in the study of self-matching properties of the
graph $\lfloor j\beta\rfloor $ against $j$ is played by the
generalized Fibonacci sequence $(G_i)_{i\in\N}$, defined
by~\eqref{ega} and~\eqref{egb}, respectively. It turns out that
shifting the argument $j$ of the function $\lfloor j\beta\rfloor$
by the integer $G_i$ results in shifting the value by $G_{i-1}$,
except of isolated mismatches with low frequency. The first
proposition is an easy consequence of the expressions of $\beta^i$
as an element of the ring $\Z[\beta]$ in the integral basis
$1,\beta$, given by~\eqref{egaa} and~\eqref{egbb}.

\begin{thm}\label{tt}
Let $\beta\in(0,1)$ satisfy $\beta^2+m\beta=1$ and let
$(G_i)_{i=0}^\infty$ be defined by~\eqref{ega}. Let $i\in\N$. Then
for $j\in\Z$ we have
$$
\bigl\lfloor \beta(j+G_i)\bigr\rfloor  = \lfloor \beta j\rfloor
+G_{i-1} \ + \ \varepsilon_{i}(j)\,,\qquad\hbox{ where }
\quad\varepsilon_{i}(j)\in\bigl\{0,(-1)^{i+1}\bigr\}\,.
$$
The frequency of integers $j$, for which the value
$\varepsilon_{i}(j)$ is non-zero, is equal to
$$
\varrho_i:=\lim_{n\to\infty}\frac{\# \{j\in\Z\mid -n\leq j\leq n,\
\varepsilon_{i}(j)\neq 0\}}{2n+1} = \beta^i\,.
$$
\end{thm}

\pf The first statement is trivial. For, we have
\begin{equation}
\begin{aligned}
\varepsilon_{i}(j)  &=\bigl\lfloor \beta(j+G_i)\bigr\rfloor  -
\lfloor \beta j\rfloor -G_{i-1} =  \bigl\lfloor \beta j-\lfloor
\beta j\rfloor + \beta G_i - G_{i-1}\bigr\rfloor  =\\[2mm]
&= \bigl\lfloor \beta j-\lfloor \beta j\rfloor +
(-1)^{i+1}\beta^i\bigr\rfloor \in \bigl\{0,(-1)^{i+1}\bigr\}\,.
\end{aligned}\label{ee}
\end{equation}
The frequency $\varrho_i$ is easily determined in the proof of
Theorem~\ref{t}.
 \pfk

In the following theorem we determine the integers $j$, for which
$\varepsilon_i(j)$ is non-zero. From this, we easily derive the
frequency of such mismatches.

\begin{thm}\label{t}
With the notation of Theorem~\ref{tt}, we have
$$
\varepsilon_i(j) = \left\{
\begin{array}{cl}
0             &\hbox{ if }\ j\notin U_i \,,\\[2mm]
(-1)^{i+1}    &\hbox{ otherwise,}
\end{array}
\right.
$$
where
$$
U_{i} = \bigl\{ kG_{i+1} + \lfloor k\beta\rfloor  G_i\,\bigm|\,
k\in\Z,\,k\neq 0\bigr\}\ \cup\
\bigl\{\tfrac{(-1)^i-1}{2}\,G_i\bigr\}\,.
$$
\end{thm}

Before starting the proof, let us mention that for $i$ even, the
set $U_i$ can be written simply as $U_{i} = \bigl\{ kG_{i+1} +
\lfloor k\beta\rfloor  G_i\,\bigm|\, k\in\Z\bigr\}$. For $i$ odd,
the element corresponding to $k=0$ is equal to $-G_i$ instead of
$0$. The distinction according to parity of $i$ is necessary here,
since unlike the paper~\cite{bunder}, we determine the values of
$\varepsilon_i(j)$ for $j\in\Z$, not only $j\geq1$.

\pf It is convenient to distinguish two cases according to the
parity of $i$.

\begin{trivlist}
\item[$\bullet$] First let $i$ be even. It is obvious
from~\eqref{ee}, that $\varepsilon_i(j)\in\{0,-1\}$ and
\begin{equation}\label{esude}
\varepsilon_i(j)=-1 \qquad\hbox{ if and only if }\qquad \beta j -
\lfloor \beta j\rfloor  \in[ 0 ,{\beta}^i )\,.
\end{equation}
Let us denote by $M$ the set of all such $j$,
$$
M=\left\{ j\in\Z \bigm| \beta j - \lfloor \beta j\rfloor  \in [ 0
,{\beta}^i ) \right\} =\left\{ j\in\Z \bigm| k + \beta j  \in
[0,\beta^i),\hbox{ for some }k\in\Z \right\}\,.
$$
Therefore $M$ is formed by the irrational parts of the elements of
the set
$$
\begin{aligned}
\left\{ k + j\beta' \bigm| k+j\beta \in [0,\beta^i) \right\} =
\Sigma[0,{\beta}^i) ={\beta'}^i \Sigma[0,1) = \\=(-\beta' G_i +
G_{i-1}) \left\{ k\beta' - \lfloor k\beta\rfloor  \bigm|
k\in\Z\right\}\,.
\end{aligned}
$$
Separating the irrational part we obtain
$$
\begin{aligned}
M&=\left\{ k G_i m  + k G_{i-1} + \lfloor k\beta\rfloor  G_{i}
\bigm| k\in\Z \right\} =\\ &=\left\{ G_i\lfloor k\beta\rfloor  + k
G_{i+1} \bigm| k\in\Z \right\} =U_i\,, \end{aligned} $$ where we
have used the equations ${\beta'}^2+m\beta'=1$ and $mG_i+G_{i-1} =
G_{i+1}$.

\item[$\bullet$] Let now $i$ be odd. Then from~\eqref{ee},
$\varepsilon_i(j)\in\{0,1\}$ and
\begin{equation}\label{eliche}
\varepsilon_i(j)=1 \qquad\hbox{ if and only if }\qquad \beta j -
\lfloor \beta j\rfloor  \in [1-\beta^i,1)\,.
\end{equation}
Let us denote by $M$ the set of all such $j$,
$$
\begin{aligned}
M&=\left\{ j\in\Z \bigm| \beta j - \lfloor \beta j\rfloor  - 1 \in
[ -{\beta}^i,0 ) \right\} =\\ &= \left\{ j\in\Z \bigm| k + \beta j
\in [-\beta^i,0),\hbox{ for some }k\in\Z \right\}\,. \end{aligned}
$$ Therefore $M$ is formed by the irrational parts of elements of
the set $$
\begin{aligned} \left\{ k + j\beta' \bigm| k+j\beta \in [-\beta^i,0)
\right\} = \Sigma[-{\beta}^i,0) ={\beta'}^i \Sigma[-1,0) =\hspace*{1cm}\\
={\beta'}^i \bigl(-1+\Sigma[0,1)\bigr)  =(\beta' G_i - G_{i-1})
\left\{ k\beta' - \lfloor k\beta\rfloor  - 1 \bigm|
k\in\Z\right\}\,.
\end{aligned} $$ Separating the irrational part we obtain $$
\begin{aligned} M&=\left\{ - k G_i m  - k G_{i-1} - \lfloor k\beta\rfloor
G_i -G_{i} \bigm| k\in\Z \right\}=\\ &= \left\{ - k G_{i+1} -
G_i\bigl(\lfloor k\beta\rfloor  + 1\bigr)  \bigm| k\in\Z
\right\}=\\ &= \left\{  k G_{i+1} + G_i\bigl(\lceil k\beta\rceil -
1\bigr) \bigm| k\in\Z \right\} =U_i\,, \end{aligned} $$ where we
have used the equation ${\beta'}^2+m\beta'=1$, $mG_i+G_{i-1} =
G_{i+1}$ and $-\lfloor -k\beta\rfloor  =\lceil k\beta\rceil$.
\end{trivlist}

Let us recall that the Weyl theorem~\cite{weyl} says that numbers
of the form $\alpha j - \lfloor \alpha j\rfloor $, $j\in\Z$, are
uniformly distributed in $(0,1)$ for every irrational $\alpha$.
Therefore the frequency of those $j\in\Z$ that satisfy $\alpha j -
\lfloor \alpha j\rfloor \in I\subset(0,1)$ is equal to the length
of the interval $I$. Therefore one can derive from~\eqref{esude}
and \eqref{eliche} that the frequency of mismatches (non-zero
values $\varepsilon_{i}(j)$) is equal to $\beta^i$, as stated by
Theorem~\ref{tt}. \pfk

If $\beta\in(0,1)$ is the quadratic unit satisfying
$\beta^2-m\beta=-1$, then the considerations are even simpler,
because the expression~\eqref{egbb} does not depend on the parity
of $i$. We state the result as the following theorem.

\begin{thm}
Let $\beta\in(0,1)$ satisfy $\beta^2-m\beta=-1$ and let
$(G_i)_{i=0}^\infty$ be defined by~\eqref{egb}. For $i\in\N$, put
$$
V_{i} = \bigl\{kG_{i+1} - (\lfloor  k\beta\rfloor +1)
G_i\,\bigm|\, k\in\Z\bigr\}\,.
$$
Then for $j\in\Z$ we have
$$
\bigl\lfloor \beta(j+G_i)\bigr\rfloor  = \lfloor \beta j\rfloor
+G_{i-1} \ + \ \varepsilon_{i}(j)\,,
$$
where
$$
\varepsilon_i(j) = \left\{
\begin{array}{cl}
0  &\hbox{ if }\ j\notin V_{i} \,,\\[2mm]
1           &\hbox{ otherwise.}
\end{array}
\right.
$$
The density of the set $U_i$ of mismatches is equal to $\beta^i$.
\end{thm}

\pf The proof follows the same lines as proofs of
Theorems~\ref{tt} and~\ref{t}. \pfk

\section{Conclusions}

One-dimensional cut-and-project sets can be constructed from
$\Z^2$ for every choice of straight lines $V_1$, $V_2$, if the
latter have irrational slopes. However, in our proof of the
self-matching properties of the Beatty sequences we strongly use
the algebraic ring structure of the set $\Z[\beta']$, and its
scaling invariance with the factor $\beta'$, namely
$\beta'\Z[\beta]=\Z[\beta']$. For that, $\beta'$ being quadratic
unit is necessary.

However, it is plausible, that even for other irrationals
$\alpha$, some self-matching property is displayed by the graph
$\lfloor j\alpha\rfloor $ against $j$. For showing that, other
methods would be necessary.


\section*{Acknowledgements}

The authors acknowledge financial support by Czech Science
Foundation GA \v{C}R 201/05/0169, by the grant LC06002 of the
Ministry of Education, Youth and Sports of the Czech Republic.



\end{document}